\documentclass[11pt,amsmath,amsthm,amsfonts,amssymb,amscd]{article}
\usepackage{epsfig}
\usepackage{relsize}
\usepackage{amsmath}
  \usepackage{paralist}
 \usepackage{psfrag, subfigure}
  \usepackage{graphics} 
  \usepackage{epsfig} 
\usepackage{graphicx}  \usepackage{epstopdf}
 \usepackage[colorlinks=true]{hyperref}
\hypersetup{urlcolor=blue, citecolor=red}

\newtheorem{theorem}{Theorem}[section]

\newtheorem{lemma}[theorem]{Lemma}

\newtheorem{definition}[theorem]{Definition}
\newtheorem{rk}{Remark}

\newcommand{\pf}{{\flushleft{\bf Proof: }}}
\newcommand{\N}{\mbox{$\mathcal{N}$}}
\newcommand{\Z}{\mbox{$Z\!\!\!Z$}}

\numberwithin{equation}{section}

\begin{document}

\title{ $C^{0}$-stability for actions implies shadowing property. }

\author{Jorge Iglesias and Aldo Portela}




\maketitle

\centerline{\scshape  Jorge Iglesias$^*$ and Aldo Portela$^*$}
\medskip
{\footnotesize
 \centerline{Universidad de La Rep\'ublica. Facultad de Ingenieria. IMERL}
   \centerline{ Julio Herrera y Reissig 565. C.P. 11300}
   \centerline{ Montevideo, Uruguay}}

\bigskip

 \centerline{(Communicated by )}

\begin{abstract}
In this paper we consider actions of the free group on two generators.
 We prove that an action of this group on a compact manifold of dimension greater or equal to two that is $C^{0}$-stable must have the shadowing property.
 We also construct both a $C^0$ and  a $C^{1}$-stable action on $S^{1}$.
\end{abstract}

\section{Introduction.}

Given a topological group $G$ and a metric space $(X,d)$,  a dynamical system is formally defined as a triplet $(G,X,\Phi )$, where $\Phi:G \times X \to X$ is a continuous function
with $\Phi (  g_1, \Phi(g_2 ,x))= \Phi(g_1 g_2 ,x)$ for all $g_1,g_2 \in G$ and for all $x \in X$. The map $\Phi$ is called an action of $G$ on $X$. It is possible to associate to
each element of $g\in G$ a homeomorphism $\Phi_g :X \to X$ where $\Phi_g(x)=\Phi(g,x)$.
For every $x \in X$ we define the orbit of $x$ as $O(x)=\{\Phi_g(x): \ g\in G\} $. A non-empty set $M \subset X$ is called minimal for the action if $\overline{O(x)}=M$ for any
$x \in M$.\\
A group $G$ is finitely generated if there exists a finite set $S\subset G$ such that for any $g\in G$ there exist $s_1,..., s_n\in S$ with $g=s_1.\cdots . s_n$.
The set $S$ is called  a finite generator of $G$. If $S$ is a finite generator of $G$ and for all $s\in S$ we have that $s^{-1}\in S$, then the set $S$ is called a finite symmetric generator.\\
For usual dynamical systems, that is when the group is  $\Z$ and the action is $\Phi (n,x)=f^{n}(x)$, we say that a sequence $\{x_n\}$ is
$\delta$-pseudotrajectory if $$d(f(x_n), x_{n+1})<\delta ,  \ \ \forall n\in Z.$$  It is possible to generalize this concept to dynamical systems $(G,X,\Phi )$ as
introduced in \cite{ot}:
A $G$-sequence in $X$ is a function $F:G\to X$. We denote this function by $\{ x_g  \} $ where $F(g)=x_g$.
Let $S$ be a finite symmetric generator of $G$. For $\delta >0$ we say that a $G$-sequence $\{ x_g  \} $  is a $\delta$-pseudotrajectory
 if
$$d(\Phi_s(x_g), x_{sg})<\delta ,\ \  \forall g\in G \mbox{ and } \forall s\in S.$$

Given a dynamical system $(G,X,\Phi)$, we say that $\Phi$  has the shadowing
property if for any $\varepsilon > 0$ there exists $\delta > 0$ such that
for any $\delta$-pseudotrajectory $\{ x_g  \} $  there exists a point $y \in X$ with
$$d(x_g, \Phi_g(y)) < \varepsilon, \ \forall g \in G.$$

A finite symmetric generator $S$ of $G$ is uniformly continuous if for every $\varepsilon >$ 0 there exists
$\delta  > 0$ such that $d(x, y) < \delta$ implies $d(\Phi_s (x),\Phi_s (y))<\varepsilon $ for every $s\in S$.
It was proved in {\cite[Proposition 1]{ot}} that the shadowing property does not depend on the finite symmetric generator $S$ if it is uniformly continuous.\\

In the case of homeomorphism and diffeomorphisms ($G=\Z$) there is a strong relationship between stability and the shadowing property. An homeomorphism $f:X\to X$ is $C^{0}$-stable
 if for every $\varepsilon >0$, there exists $\delta > 0$ such that
if $g$ is another homeomorphism on $X$ with $d_{C^{0}}(f, g):=sup_{x\in X}\{ d(f(x),g(x)\} < \delta$ then there exists a continuous and surjective
function $h : X \to X $ with $hg=fh $  and $d_{C^{0}}(h,id)<\varepsilon.$ When $X$ is a manifold and $f:X\to X$ is a diffeomorphism, we say that $f$ is $C^{1}$-stable if there exists $\delta > 0$ such that $g$ is another diffeomorphism with
$d_{C^{1}}(f, g):= Max\{ d_{C^{0}}(f,g), d_{C^{0}}(Df,Dg)\}< \delta$, then there exists a homeomorphisms $h : M \to M $ with $hg=fh$.
In \cite{w} was proved that every $C^{0}$-stable homeomorphism in a compact manifold of dimension greater or equal to two, has the shadowing property. This property is also valid for $C^{1}$-stable diffeomorphisms.
The concepts of $C^{0}$ and $C^{1}$ stability can by generalized for actions as follows:
let $S$ be a finite generating set for a
group $G$. We denote by $Act(G, X)$ the set of actions in $X$, and we define a metric $d^{0}_S$ on $Act(G, X)$ by
$$d^{0}_S(\widetilde{\Phi}, \Phi) := Max_{  s\in S} \{  d_{C^{0}}(\widetilde{\Phi}_s, \Phi_s)\}, \mbox{ for } \widetilde{\Phi}, \Phi \in  Act(G, X).$$
 We say that an action $\Phi \in Act(G, X)$ is $C^{0}$- stable if for every $\varepsilon >0$, there exists $\delta > 0$ such that
if $\widetilde{\Phi}\in  Act(G, X)$  with $d^{0}_S(\widetilde{\Phi}, \Phi) < \delta$ then there exists a continuous and surjective
function $h : X \to X $, with $h\widetilde{\Phi}_g=\Phi_gh ,$ for every $g \in G$ and $d_{C^{0}}(h,id)<\varepsilon.$
It is not hard to prove that this definition does not depend on the generator $S$ if $X$ is a compact set (see \cite[Lemma 2.2]{cl}).\\
When $X$ is a manifold $M$, let  $C^{1}(M)$ be the set of $C^{1}$ maps $f:M\to M$. Let $S$ be a finite generator of $G$.
 Let $$Act^{1}(G, M)=\{ \Phi\in  Act(G, M): \ \Phi_s \in C^{1}(M) \mbox{ for all }  s\in S  \}. $$
and consider
$$d^{1}_S(\widetilde{\Phi}, \Phi) := Max_{  s\in S} \{  d_{C^{1}}(\widetilde{\Phi}_s, \Phi_s)\}, \mbox{ for } \widetilde{\Phi}, \Phi \in  Act^{1}(G, X).$$

An action $\Phi\in Act^{1}(G, M)$ is $C^{1}$- stable if there exists $\delta > 0$ such that
if $\widetilde{\Phi} \in Act^{1}(G, M)$ with
$d^{1}_{S}(\widetilde{\Phi}, \Phi) < \delta$, then there exists an homeomorphism $h : M \to M $ with $h\widetilde{\Phi}_g=\Phi_gh ,$ for every $g \in G$.\\
We denote by $F_2$ the free group generated by two elements. Let us state our main result:\\
{\bf{Theorem A}}
{\it{ Let $(F_2,M,\Phi)$ be a dynamical system with $M$ a compact manifold of dimension greater or equal to two. If $\Phi$ is $C^{0}$-stable then $\Phi$ has the shadowing property.}}\\

An action $\Phi \in Act(G, X)$ is expansive if there exist $\alpha >0$ such that for every $x,y\in X$ with $x\neq y$ there exists $g\in G$ such that $d(\Phi_g(x),\Phi_g(y))>\alpha$.
Let $X$ be a  compact metric space. In \cite{cl} was proved that if the action $\Phi \in Act(G, X)$ is expansive and has the shadowing property then $\Phi$ is $C^{0}$-stable. In
section 2,  we construct a $C^{0}$ and $C^{1}$-stable action of  $G=F_2$ on $S^{1}$ that  is not expansive. This shows that the expansive property is not a necessary condition for the stability.\\

\section{Example of a $C^{0}$ and $C^{1}$-stable action. }
In this section, we consider the free group $F_2$ with finite symmetric generator $S=\{a,a^{-1},b,b^{-1}\}$.
We are going to construct a $C^{0}$ and $C^{1}$-stable action $\Phi$ in $\mathcal{S}^2 $ whose minimal set $K$ is a Cantor set. The generator of the action will be $\Phi_a $ and $\Phi_b$ where $\Phi_a, \Phi_b :\mathcal{S}^2 \to \mathcal{S}^2 $ will be defined later.\\

\begin{figure}[ht]
\begin{center}
\caption{\label{phia2}}
\psfrag{phia}{$f_a$}
\psfrag{ff}{$f^{-1}_a$}
\psfrag{sa}{$S_a$}
\psfrag{ib}{$I_b$}
\psfrag{ibb}{$I_{b^{-1}}$}
\psfrag{ia}{$J_a$}
\psfrag{iia}{$J_{a^{-1}}$}
\psfrag{ii}{$+\infty$}
\psfrag{0}{$0$}
\includegraphics[scale=.2]{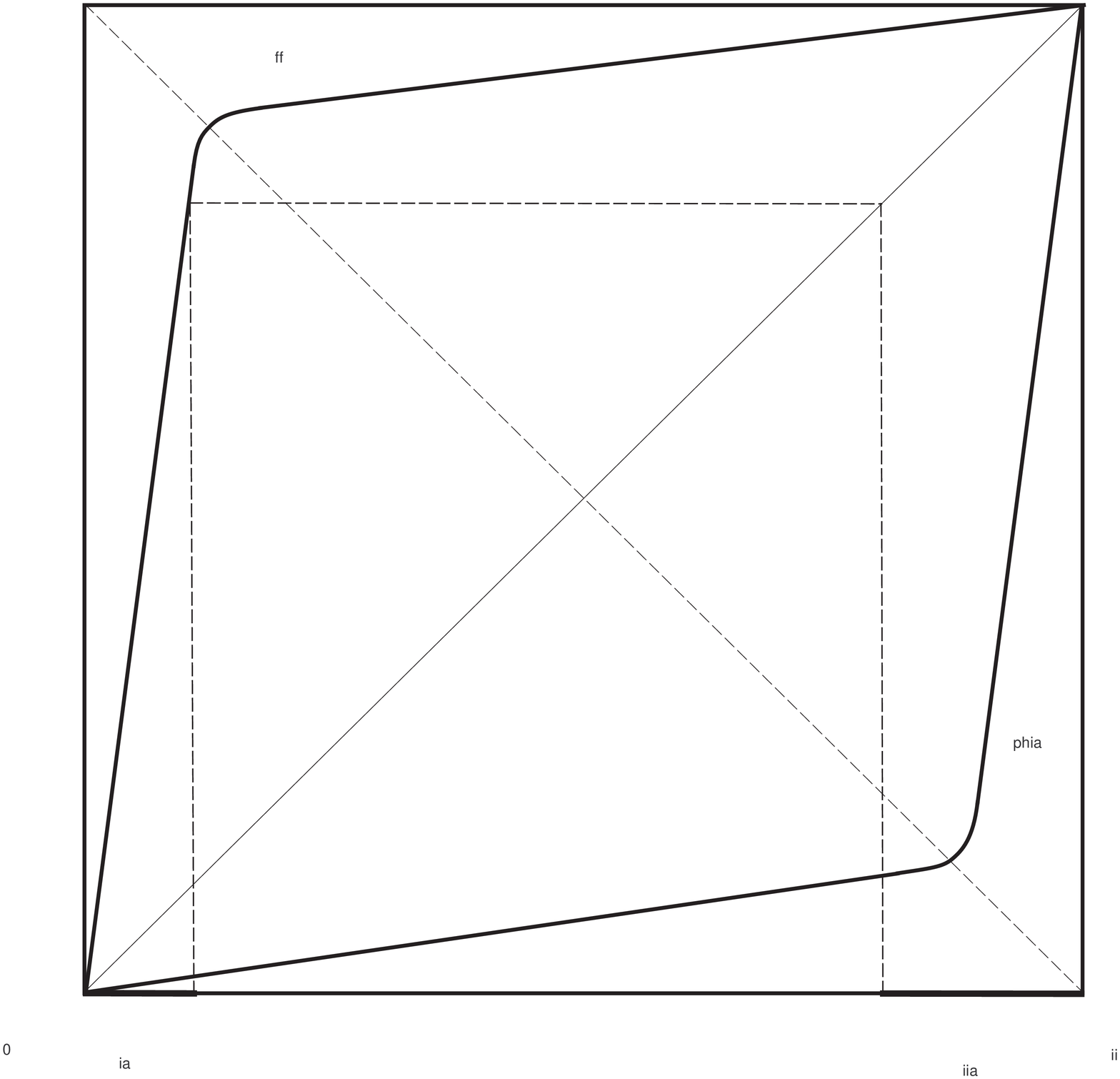}
\end{center}
\end{figure}

We begin defining a function $f_a:[0,+\infty]\to [0,+\infty]$, $f_a\in C^{1}$, as figure \ref{phia2}. Considerer the intervals $J_a$ y  $J_{a^{-1}}= \overline{f^{-1}_a(J_a)}^{c}$ such that (see figure \ref{phia2}):

\begin{enumerate}
\item  $|f^{-1}_a(x)-f^{-1}_a(y)|>2| x-y   |$ for all $x,y\in J_a$ (and consequently $|f_a(x)-f_a(y)|<1/2| x-y   |$ for all $x,y\in J_a$ ).

\item   $|f^{-1}_a(x)-f^{-1}_a(y)|<1/2| x-y   |$ for all $x,y\in J_{a^{-1}}$.

\end{enumerate}

Now we defining $\Phi_a :\mathcal{S}^2 \to \mathcal{S}^2 $ (in polar coordinates) and the intervals $I_a$ y $I_{a^{-1}} $ such that
\begin{itemize}
\item $ \Phi_a (r,\theta )=(f_a(r), \theta).    $
\item $I_a=\{ (r,\theta): \ \ r\in J_a \mbox { and } \theta\in [0,2\pi )   \}$ and $I_{a^{-1}}=\{ (r,\theta): \ \ r\in J_{a^{-1}} \mbox { and } \theta\in [0,2\pi )   \}$
\end{itemize}

Note that  $I_{a^{-1}}= \overline{\Phi^{-1}_a(I_a)}^{c}$.\\

It is clear that $\Phi_a:\mathcal{S}^2 \to \mathcal{S}^2 $ is a north-south pole diffeomorphism, with $\Omega (\Phi_a)=\{ N_a,S_a  \}$,  ( $S_a\in I_a$, $N_a\in I_{a^{-1}}$  ) and verifies
\begin{itemize}
\item  $||\Phi^{-1}_a(x)-\Phi^{-1}_a(y)||>2|| x-y   ||$ for all $x,y\in I_a$ (and consequently $||\Phi_a(x)-\Phi_a(y)||<1/2|| x-y   ||$ for all $x,y\in I_a$ ).

\item   $||\Phi^{-1}_a(x)-\Phi^{-1}_a(y)||<1/2|| x-y   ||$ for all $x,y\in I_{a^{-1}}$.
\end{itemize}

Now we consider  $\Phi_b:\mathcal{S}^2 \to \mathcal{S}^2 $ north-south pole diffeomorphism, with $\Omega (\Phi_b)=\{ N_b,S_b  \}$ defined analogously to $\Phi_a$.
Note that it is possible to construct $\Phi_a$ and $\Phi_b$ so that they verify the following properties:

\begin{enumerate}

\item $ \overline{I}_a\cap \overline{I}_s=\emptyset  $, for $s\in \{ b,b^{-1}\}$ and $ \overline{I}_{a^{-1}}\cap \overline{I}_s=\emptyset  $, for $s\in \{ b,b^{-1}\}$

\item $||\Phi_a(x)-\Phi_a(y)||<1/2|| x-y   ||$ for all $x,y\in  I_b\cup I_{b^{-1}} $.

\item $||\Phi^{-1}_a(x)-\Phi^{-1}_a(y)||<1/2|| x-y   ||$ for all $x,y\in  I_b\cup I_{b^{-1}} $.

\item $||\Phi_b(x)-\Phi_b(y)||<1/2|| x-y   ||$ for all $x,y\in  I_a\cup I_{a^{-1}} $.

\item $||\Phi^{-1}_b(x)-\Phi^{-1}_b(y)||<1/2|| x-y   ||$ for all $x,y\in  I_a\cup I_{a^{-1}} $.

\end{enumerate}

\begin{figure}[h]
\psfrag{ia}{$I_{a}$}

\psfrag{ib}{$I_{b}$}

\psfrag{iaa}{$I_{a^{-1}}$}

\psfrag{ibb}{$I_{b^{-1}}$}

\psfrag{aaiaa}{$\Phi_a^{-1} ( I_{a^{-1}}) $}

\psfrag{aaib}{$\Phi_a^{-1} ( I_{b}) $}

\psfrag{aaibb}{$\Phi_a^{-1} ( I_{b^{-1}}) $}

\psfrag{bbibb}{$\Phi_b^{-1} ( I_{b^{-1}}) $}

\psfrag{bbiaa}{$\Phi_b^{-1} ( I_{a^{-1}}) $}

\psfrag{bbia}{$\Phi_b^{-1} ( I_{a}) $}

\psfrag{bib}{$\Phi_b ( I_{b}) $}

\psfrag{biaa}{$\Phi_b ( I_{a^{-1}}) $}

\psfrag{bia}{$\Phi_b ( I_{a}) $}

\psfrag{aia}{$\Phi_a( I_{a}) $}

\psfrag{aib}{$\Phi_a ( I_{b}) $}

\psfrag{aibb}{$\Phi_a ( I_{b^{-1}}) $}

\begin{center}
\caption{\label{figura44}}
\subfigure[]{\includegraphics[scale=0.27]{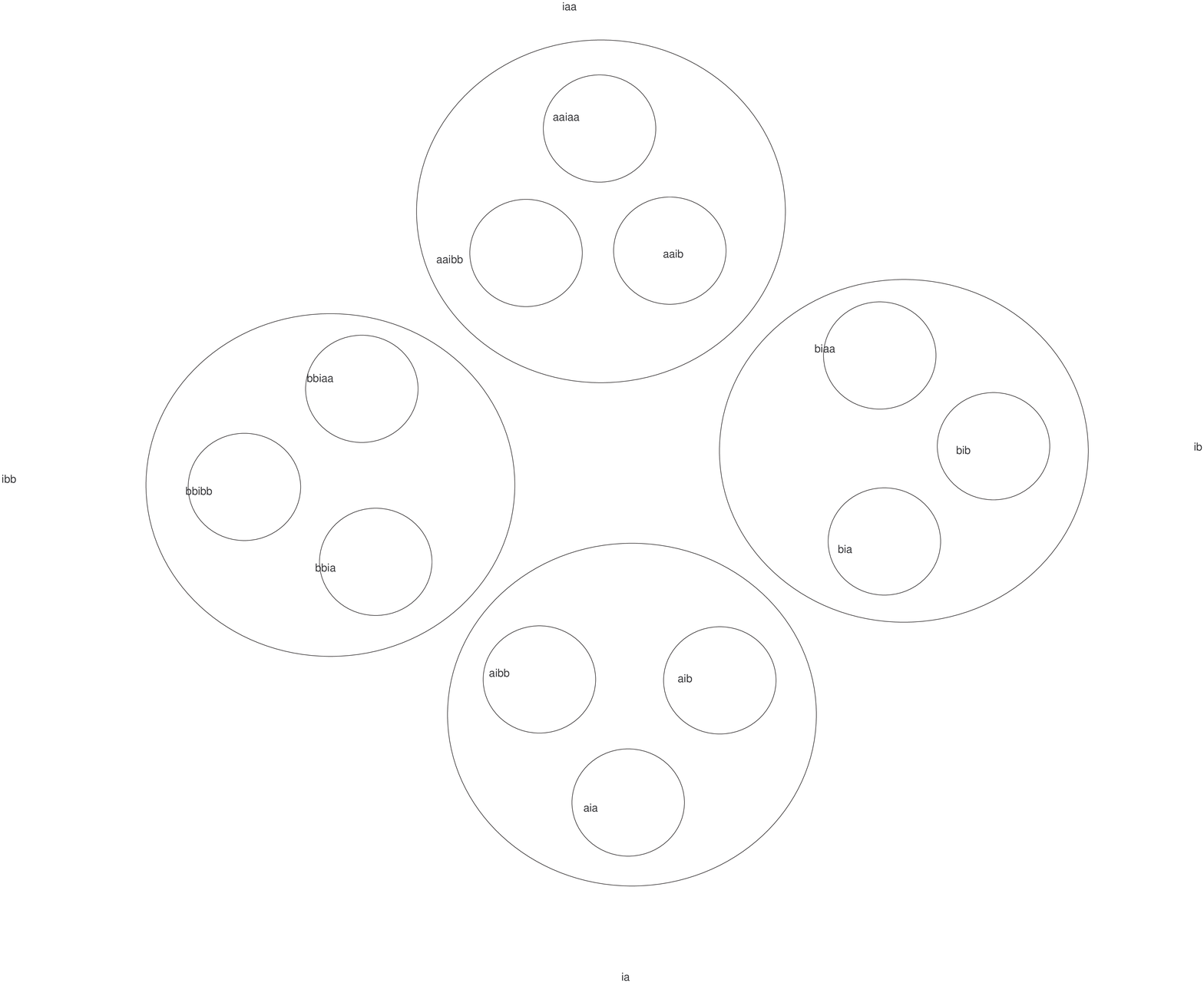}}
\end{center}
\end{figure}

Let $\Phi$ be the action generated for $\Phi_a$ and $\Phi_b$. The figure \ref{figura44} shows the dynamic of the action of $\Phi $ on the intervals $I_a, I_{a^{-1}}, I_b$ and $I_{b^{-1}}$.

 The following properties are very useful for our purpose. Since they are not hard to prove we omit its proof.

\begin{rk}\label{rk1}
\begin{enumerate}
\noindent\item If $s,s^{'}\in \{ a,a^{-1},b,b^{-1}\}$ then $$||\Phi_s(x)-\Phi_s(y)||<1/2|| x-y   || \ \forall x,y\in I_{s^{'}} \mbox{ with } s^{'}\neq s^{-1}.$$
\item If $s\neq s{^{'}{^{-1}}}$ then $\Phi_s (I_{s^{'}})\subset I_s$.
\item $\Phi_s(I_s)\cap I_{s^{-1}}=\emptyset$ for all $s\in \{a,b,a^{-1},b^{-1}\}$.\\
\end{enumerate}
\end{rk}

Let $\{A_n\}$ be such that
 $$A_0=I_a\cup I_{a^{-1}}\cup I_b\cup I_{b^{-1}}  \mbox{ and }$$ $$  A_{n+1}= \left[ \Phi_a(A_n)\cap \Phi_{a^{-1}}(A_n)\cap \Phi_b(A_n)\cap \Phi_{b^{-1}}(A_n)   \right] \cap A_n.$$
 Note that

 \begin{itemize}
\item For any $n\in\N$, $A_n$ has $4.3^{n}$ connected components and $\overline{A_{n+1}}\subset int(A_n)$.
\item The diameter of the connected components fo $A_n$ goes to zero when $n$ goes to infinity.
\item The Cantor set  $K=\bigcap_{n\geq 1}A_n$ is a minimal set for the action $\Phi$ generated for $\Phi_a$ y $\Phi_b$ (see \cite{n}).
\end{itemize}

{\bf{Construction of the semiconjugacy.}}

For the action $\Phi$ we can think that the set $K$ is an "attractor set" and that $A_0^{c}$ plays the role of a "fundamental domain".
To prove that $\Phi$ is a $C^{0}$-stable action we will take $\delta>0$ such that if $d^{0}_S(\widetilde{\Phi}, \Phi  )<\delta$ then the dynamic of $\widetilde{\Phi}$ is given as figure  \ref{figura44}.
So we define $h:\widetilde{A}^{c}_0\to A^{c}_0$ close to identity, we extend it dynamically to $\mathcal{S}^2\setminus K$ and applying Lemma \ref{extension_lemma} it can be extended to $\mathcal{S}^2$.

Fix the discs $I_a$ and $I_b$.
Let $\alpha= min\{ d_H ( \overline{I_s}, \overline{I_{s^{'}}}  ), \  s,s^{'} \in\{a,a^{-1}, b,b^{-1}\}   \mbox{ and } s\neq s^{'}  \}$ be, where  $d_H$ is the Hausdorff distance. Note that if $x,y\in K$, $x\neq y$, then there exists $g\in G$ such that $d(\Phi_g(x),\Phi_g(y))>\alpha $.
For any other action $\widetilde{\Phi}$ considerer  $\widetilde{I}_a= I_a$,  $\widetilde{I}_b=I_b$,   $ \widetilde{I}_{a^{-1}}= \overline{\widetilde{\Phi}^{-1}_a( \widetilde{I}_a)}^{c}$ and  $\widetilde{I}_{b^{-1}}= \overline{\widetilde{\Phi}^{-1}_b(\widetilde{I}_b)}^{c}$. Let $$\widetilde{A}_0= \widetilde{I}_a \cup \widetilde{I}_{a^{-1}}   \cup \widetilde{I}_b  \cup \widetilde{I}_{b^{-1}},$$ $$\widetilde{A}_{n+1}= \left[ \widetilde{\Phi}_a( \widetilde{A}_n)\cap \widetilde{\Phi}_{a^{-1}}(\widetilde{A}_n)\cap \widetilde{\Phi}_b(\widetilde{A}_n)\cap \widetilde{\Phi}_{b^{-1}}(\widetilde{A}_n)   \right] \cap \widetilde{A}_n$$
$$\mbox{ and we call } K^{c}_{\widetilde{\Phi}}=\cap \widetilde{A}_n. $$

\begin{definition}\label{definicion_delta}{\bf{Neighbourhood of stability.}}

Given $\varepsilon>0$, $\varepsilon <\alpha /8$, let $\delta >0$ be, $\delta <\varepsilon $, such that if $\widetilde{\Phi}$ is an action with $ d^{0}_S(\Phi, \widetilde{\Phi}  )<\delta $, then:

\begin{enumerate}\label{def_phi}

\item $\overline{ \widetilde{I}_s}\cap \overline{\widetilde{I}_{s^{'}}}= \emptyset $ for $s,s^{'}\in \{ a,a^{-1},b,b^{-1}\}$ with $s\neq s^{'}$.

\item If $s\neq s{^{'}{^{-1}}}$ then $\overline{\widetilde{\Phi}_s (\widetilde{I}_{s^{'}}) }\subset \widetilde{I}_s$.

\item If $I$ is a connected  component of $K_{\widetilde{\Phi}}$, then $diam(I)<\varepsilon /2.$

\item For all $s\in\{a,b,a^{-1},b^{-1}\}$, $d( \widetilde{\Phi}^{-1}_{s} \Phi_{s}(x),x )<\varepsilon/2$ and $d( \Phi^{-1}_{s} \widetilde{\Phi}_{s}(x),x )<\varepsilon/2$ .

\item For all $x,y\in  I_{s^{'}}\cup \widetilde{I}_{s^{'}} $, $s\neq (s^{'})^{-1}$,  $||\Phi_s(x)-\Phi_s(y)||<1/2|| x-y   ||$ .

\end{enumerate}

\end{definition}

\vspace{1cm}

\begin{figure}[ht]
\begin{center}
\caption{\label{perturbado}}
\psfrag{phia}{$\widetilde{f}_a$}
\psfrag{i}{$J$}

\includegraphics[scale=.12]{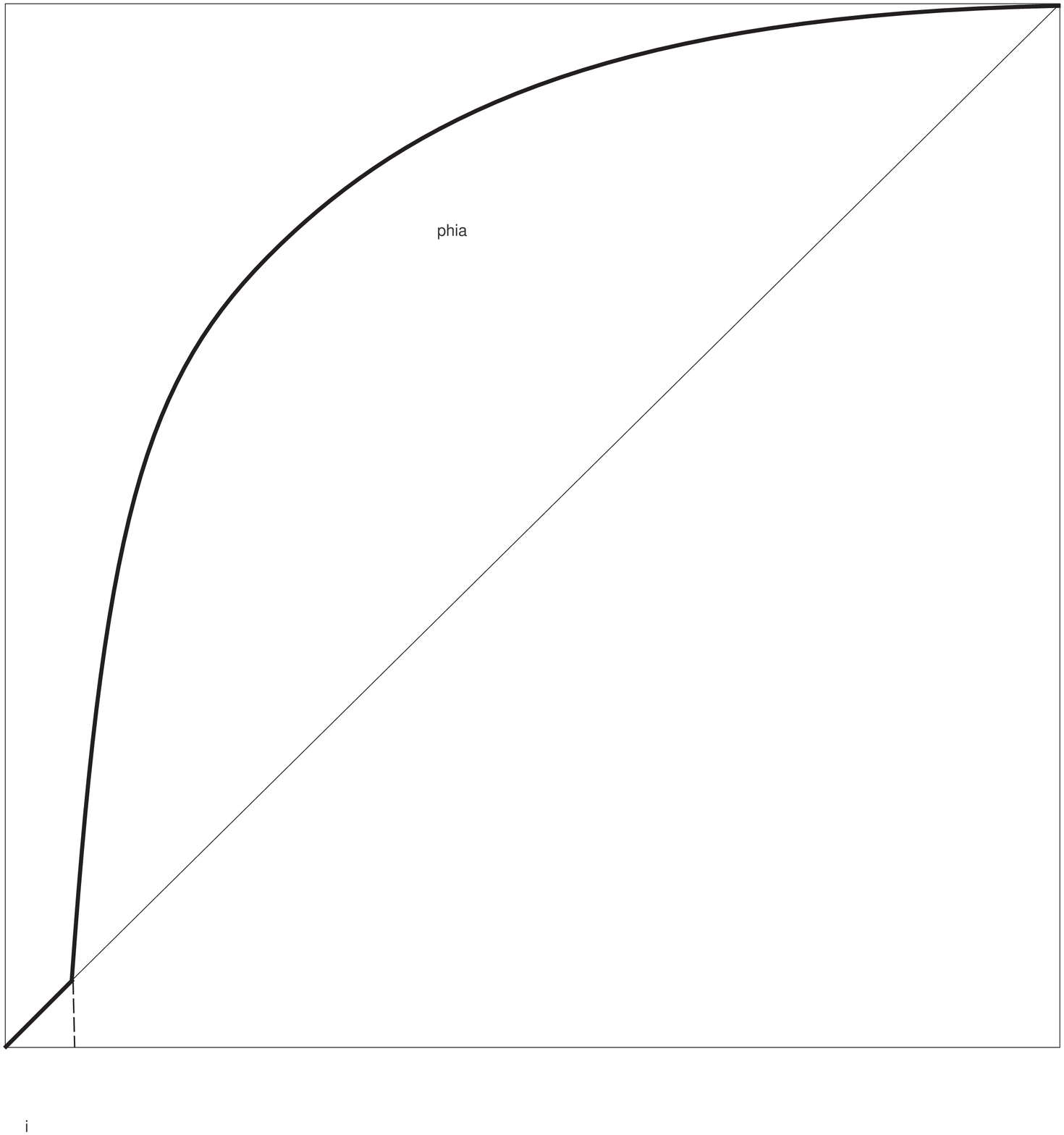}
\end{center}
\end{figure}

Note that if $\widetilde{\Phi}_a (r,\theta)= (\widetilde{f}_a(r),\theta)$ is a $C^{0}$-perturbed of $\Phi_a$, where $\widetilde{f}_a$ is as figure \ref{perturbado}, then $K_{\widetilde{\Phi}}$ is not a Cantor set.

For the following lemmas, the action  $\widetilde{\Phi}$ satisfies the five properties above.

The following remark will be very helpful.

\begin{rk}\label{rk12}
\begin{enumerate}
\item $\widetilde{\Phi}_s(\widetilde{I}_s)\cap \widetilde{I}_{s^{-1}}=\emptyset$ for all $s\in \{a,b,a^{-1},b^{-1}\}$.\\
\end{enumerate}
\end{rk}

Some of our proofs are by induction in the length of the elements $ g\in G $. Thus we need to define the length of an element $g\in G $.
The elements of length one are $ a,a^{-1},b$ and $b^{-1}$. The elements of length $n$ are obtained from the elements of length $n-1$ as follows: Let $g=s_{n-1}....s_{1}$ be an element of lengths $n-1$ with $s_j\in \{ a,a^{-1},b,b^{-1}\}$. Then the elements of length $n$ generated by $g$ are $g^{'}=s.g$ with $s\neq  (s_{n-1})^{-1}$. The length of $g$ is denoted by $|g|$.
It is clear that an element $g \in G $ can be written from $ S $ in different ways, for example $ g = gaa ^ {- 1} $.
Note that if $g=s_{{n}}....s_{{1}}$ with $s_j\in \{ a,a^{-1},b,b^{-1}\}$, then  $|g|\leq n$. We say that $g=s_{{n}}....s_{{1}}$ is written in its normal form if $|g|= n$. It is easy to prove that the normal representation is unique. From now on we will consider $g\in G$ written in its normal form.\\

\begin{lemma}\label{lema_importante}
Let $\widetilde{J}$ be a connected component of $\widetilde{A}_n\setminus \widetilde{A}_{n+1}$ with $\widetilde{J}\subset \widetilde{I}_s$ and $g\in F_2$, $g=s_j...s_1$ with $s_1\neq s^{-1}$. Then $\widetilde{\Phi}_g(\widetilde{J})\subset \widetilde{A}_{n+j}\setminus \widetilde{A}_{n+1+j}$.
\end{lemma}
\pf
 As $\widetilde{J}\subset \widetilde{I}_s$ and  $s_1\neq s^{-1}$, by  item 2) in the election of $\delta$, $\widetilde{\Phi}_{s_{1}}(\widetilde{J})\subset \widetilde{I}_{s_{1}}. $
   On the other hand, as $\widetilde{J}$ is a connected component of $\widetilde{A}_n\setminus \widetilde{A}_{n+1}$   then $\widetilde{\Phi}_{s_{1}}(\widetilde{J}) $  is a connected component of $\widetilde{A}_{n+1}\setminus \widetilde{A}_{n+2}$. Again as $s_2\neq s_1^{-1}$, because $g$ is written in its normal form,
then $\widetilde{\Phi}_{s_{2}} \widetilde{\Phi}_{s_{1}}(\widetilde{J}) \subset   \widetilde{I}_{s_{2}}$. As $\widetilde{\Phi}_{s_{1}}(\widetilde{J}) $  is a connected component of $\widetilde{A}_{n+1}\setminus \widetilde{A}_{n+2}$ then $\widetilde{\Phi}_{s_{2}} \widetilde{\Phi}_{s_{1}}(\widetilde{J})$   is a connected component of $\widetilde{A}_{n+2}\setminus \widetilde{A}_{n+3}$.   Reasoning inductively we obtain the thesis.\\

\begin{rk}\label{lemaJ}
\begin{enumerate}
\item  Given $\widetilde{J}$ a connected component of $\widetilde{A}_n\setminus \widetilde{A}_{n+1}$, there exists $g\in F_2$ such that $\widetilde{\Phi}_g(\widetilde{J})\subset \widetilde{A}_0^{c}$.

\item Given $\widetilde{J}$ a connected component of $\widetilde{A}_n\setminus \widetilde{A}_{n+1}$ with $\widetilde{J}\subset \widetilde{I}_s$. If $s^{'}\neq s^{-1}$ then $\widetilde{\Phi}_{s^{'}}(\widetilde{J})\subset \widetilde{A}_{n+1}$.
\end{enumerate}
\end{rk}

\begin{lemma}\label{lemma3}
  Let $\widetilde{J}$ be a connected component of $\widetilde{A}_n\setminus \widetilde{A}_{n+1}$,  with $n\geq 0$. Then there is a unique element $g_{_{\widetilde{J}}}$ of minimum length such that $ \widetilde{\Phi}_{g_{_{\widetilde{J}}}}(\widetilde{J})\subset  \widetilde{A}_0^{c}$.
\end{lemma}

\pf
By Remark \ref{lemaJ} item 1., there exists  $g_{_{\widetilde{J}}}$ of minimum length such that  $ \widetilde{\Phi}_{g_{_{\widetilde{J}}}}(\widetilde{J})\subset\widetilde{A}_0^{c}$. Let $g_{_{\widetilde{J}}}=s_r...s_1$. Let $s$ be such that $\widetilde{J}\subset \widetilde{I}_s$. By Lemma \ref{lema_importante}  $s_1= s^{-1}$.  If $n=0$ then  $\widetilde{\Phi}_{g_{_{\widetilde{J}}}}(\widetilde{J})\subset  \widetilde{A}_0^{c}$ and therefore  $g_{_{\widetilde{J}}}=s^{-1}$.

 If $n> 0$ then $ \widetilde{\Phi}_{s^{-1}}(\widetilde{J})$ is a connected component of  $\widetilde{A}_{n-1}\setminus \widetilde{A}_{n} $.
 Reasoning as in the previous case, $ s_2 $ is determined. This proves the uniqueness of $g_{_{\widetilde{J}}}$. \\

Next, we will prove that for every action  $\widetilde{\Phi}\ $ ( that verifies the five conditions above) there exists a semiconjugacy $h:\mathcal{S}^2 \to \mathcal{S}^2 $ such that $h\widetilde{\Phi}_g=\Phi_g h$ for all $g\in F_2$.\\

\begin{rk}\label{definicion_h}

 As $ \partial \widetilde{A}_0^{c}= \partial \widetilde{I}_a \cup    \partial \widetilde{I}_{a^{-1}}  \cup  \partial \widetilde{I}_b  \cup \partial \widetilde{I}_{b^{-1}} $, therefore

\begin{itemize}

\item If  $ \Phi_{g}  (  \widetilde{A}_0^{c} ) =  \widetilde{A}_{0}^{c}$ then $g=e$.

\item Let $x\in \partial \widetilde{A}_0^{c}$ be with  $\Phi_g(x)\in \partial \widetilde{A}_0^{c}$, and $g\neq e$.
If $x\in \partial \widetilde{I}_s$ then $g=s^{-1}$ and $\Phi_g(x)\in \partial  \widetilde{I}_{s^{-1}}$.

\end{itemize}

\end{rk}

We need that the following property to be satisfied. For this, we reduce $\delta$ as necessary so that in addition to complying with properties list in definiton \ref{definicion_delta}, it verify the following additional property:\\

Given $\varepsilon>0$ (as in Definition \ref{definicion_delta}), there exists $h:\overline{\widetilde{A}^{c}_0} \to \overline{A^{c}_0} $ homeomorphisms with $d(h(x),x)<\varepsilon $ and if $x\in \partial \widetilde{I}_s$ and
 $\widetilde{\Phi}_{s^{-1}} (x)\in \partial \widetilde{I}_{s^{-1}}$ then $h(x)=\Phi^{-1}_{s^{-1}}h\widetilde{\Phi}_s(x)$. This is possible by Remark \ref{definicion_h}.

By Lemma \ref{lemma3} we can extend $h$ dynamically to any connected component of $K^{c}_{\widetilde{\Phi}}$.

\begin{lemma}
 For all $x\in K^{c}_{\widetilde{\Phi}}$, $d(h(x),x)<\varepsilon$.
\end{lemma}
\pf

Let $x\in K^{c}_{\widetilde{\Phi}}$. If $x\in \widetilde{A}^{c}_0$ then, by definition of $h$, $d(h(x),x)<\varepsilon$.

  If $x\in K^{c}_{\widetilde{\Phi}}   \setminus \widetilde{A}^{c}_0$, then there exists $n\in \N$ such that $x\in\widetilde{A}_n\setminus \widetilde{A}_{n+1}$. The proof will be done by induction in $n$. If $n=0$, $x$  belong   $\widetilde{J}$,  connected component of $\widetilde{A}_0\setminus \widetilde{A}_1$ with $\widetilde{J}\subset I_{s}$ for some $s\in\{a,a^{-1},b,b^{-1}\}$.
 Then   $\widetilde{\Phi}_{s^{-1}}(\widetilde{J})\subset (\widetilde{A}_0)^{c}$.
As $h$ is dynamically defined,  so $h(x) = \Phi_sh \widetilde{\Phi}_{s^{-1}}(x)$. As $\widetilde{\Phi}_{s^{-1}}(x)\in \widetilde{A}^{c}_0 $, then $d( h\widetilde{\Phi}_{s^{-1}}(x),  \widetilde{\Phi}_{s^{-1}}(x)  )<\varepsilon$.

By item 5 of definition \ref{definicion_delta},  $d(     \Phi_s h\widetilde{\Phi}_{s^{-1}}(x), \Phi_s \widetilde{\Phi}_{s^{-1}}(x)  )<\varepsilon /2$. By item 4 of definition \ref{definicion_delta}, we have that $d( \Phi_s \widetilde{\Phi}_{s^{-1}}(x) ,x  )<\varepsilon/2$ then $d(     \Phi_s h\widetilde{\Phi}_{s^{-1}}(x), x )<\varepsilon $, therefore $d(h(x),x)<\varepsilon$.

Reasoning inductively, we finish the proof.

The following is a general lemma. As an application, $h$ is extended to $ \mathcal {S}^2 $.

\begin{lemma}\label{extension_lemma}{\bf{ Extension lemma}}.\\
Let $(X,d)$ be a compact metric space and $\Phi$, $\widetilde{\Phi}$ two actions of a group $G$. Let $X_1\subset X$ and $X_2\subset X$ be such that
 $\Phi_g(X_1)=X_1$ and  $\widetilde{\Phi}_g(X_2)=X_2$ for all $g\in G$. There exists $\alpha >0$ such that:
\begin{enumerate}
 \item For all $x,y\in \overline{X_1}\setminus X_1$, $x\neq y$, there exists $g\in G$ such that $d(\Phi_g(x),\Phi_g(y)  )>\alpha$.\\
\item Given $\varepsilon >0$, $\varepsilon <\alpha /8$, there exists a homeomorphism $h:X_2\to X_1$,  $h\widetilde{\Phi}_g=\Phi_gh$ for all $g\in G$, with $d(h(x),x)<\varepsilon$.
\end{enumerate}

Then, there exists a semiconjugacy $\overline{h}:\overline{X}_2\to \overline{X}_1$, $\overline{h}\widetilde{\Phi}_g= \Phi_g \overline{h}$ for all $g\in G$, with $d(\overline{h}(x),x)\leq \varepsilon $ and $\overline{h}|_{X_{2}}=h$.\\

Moreover, if for all $x,y\in \overline{X_2}\setminus X_2$, $x\neq y$, there exists $g\in G$ such that $d(\widetilde{\Phi}_g(x),\widetilde{\Phi}_g(y)  )>\alpha$, then  $\overline{h}:\overline{X}_2\to \overline{X}_1$ is a homeomorphisms.

\end{lemma}

\pf

  Let $x\in \overline{X_2}\setminus X_2$  and $\{x_n\}$ with $x_n\to x$ and $x_n\in X_2$ for all $n\in \N$. Suppose that $\{h(x_n)\}$ accumulate in points  $y,z$ with $y\neq z$. Let  $\{y_n\}$,  $\{z_n\}$ be with $y_n, z_n\to x$ and  $h(y_n)\to y$,  $h(z_n)\to z$. Note that $y,z\in \overline{X_1}\setminus X_1$ because $h$ is a homeomorphisms.  As $y\neq z$ there exists $\Phi_g$ such that $d( \Phi_g (y),\Phi_g(z))>\alpha $. So, for $n$ big enough you have  $d( \Phi_g (h(y_n)),\Phi_g(h(z_n)))>\alpha $. As $\Phi_g h(y_n)=h \widetilde{\Phi}_g (y_n)$, $\Phi_g h(z_n)=h \widetilde{\Phi}_g (z_n)$ then $d( h \widetilde{\Phi}_g (y_n),h \widetilde{\Phi}_g (z_n))>\alpha $. Since $d(h(x),x)<\varepsilon < \alpha/8$, we have that  $$ d( \widetilde{\Phi}_g (y_n),\widetilde{\Phi}_g (z_n))>\alpha-2\varepsilon >\frac{3\alpha}{4} \cdot $$
This is a contradiction because $y_n, z_n\to x$.\\

So, we define  $\overline{h}(x)=\lim h(x_n)$. It is easy to prove that $\overline{h}$ is continuous, surjective, semiconjugacy ($\overline{h}\widetilde{\Phi}_g=\Phi_g\overline{h}$) and $d(\overline{h}(x),x)\leq \varepsilon$.
This proves the first part of the lemma.\\

 Let's prove the second parte of lemma.\\ For this it is enough to prove that $\overline{h}$ is injective. Suppose that $\overline{h}(x)=\overline{h}(y)$ with $x\neq y$. Note that  $x,y$ must belong to $\overline{X_2}\setminus X_2$.\\

As $x\neq y$ there exists $\widetilde{\Phi}_g$ such that $d( \widetilde{\Phi}_g (x),\widetilde{\Phi}_g(y))>\alpha $. As $d(\overline{h}(x),x)\leq \varepsilon $ then
$d( \overline{h}\widetilde{\Phi}_g (x),\overline{h}\widetilde{\Phi}_g(y))\geq \frac{3\alpha}{4}$.  Then $d( \Phi_g \overline{h}(x),\Phi_g \overline{h}(y)) \geq \frac{3\alpha}{4} $. This is a contradiction because  $\overline{h}(x)=\overline{h}(y)$. Therefore $\overline{h}$ es injective. So $\overline{h}$ is a homeomorphism.\\

Now, by the Lemma \ref{extension_lemma} we can extend $h$ to $\partial K^{c}_{\widetilde{\Phi}}$.

  To finish, it is necessary to extend $h$ to $K_{\widetilde{\Phi}}$. Let $I$ be a connected component of $K_{\widetilde{\Phi}}$. Let $x,y\in\partial I\subset  \partial K_{\widetilde{\Phi}}^c$.
Note that the map $h$ is defined in the points $x,y$.
   We will prove that $h(x)=h(y)$. If $h(x)\neq h(y)$, as $h(x)$ and $ h(y)$ are in $K$, then there exists $g\in G$ such that $d(\Phi_g(h(x)),\Phi_g(h(y)))>\alpha $. Therefore

$d(h \widetilde{\Phi}_g(x),h\widetilde{\Phi}_g(y))>\alpha $. As $d(h(x),x)\leq \varepsilon$ then $d( \widetilde{\Phi}_g(x),\widetilde{\Phi}_g(y))>3\alpha /4 $.
As  $\widetilde{\Phi}_g(x)$ and $\widetilde{\Phi}_g(y)$ belong to the same connected component of $K_{\widetilde{\Phi}}$, by
 item 4. from  definition of $\delta$, we obtain a contradition.

So $h$ is constant in $\partial I$ and we define $h(I)=h(x)$. 

\vspace{.5cm}

\section{$C^{1}$-stability of action $\Phi$.}

We considerer  $Act^{1} (F_2,\mathcal {S}^{2})$ the set of  $C^{1}$ action with the distance
 $$d^{1}_S(\widetilde{\Phi}, \Phi) := Max_{ _{  s\in S}} \{  d_{C^{1}}(\widetilde{\Phi}_s, \Phi_s), \mbox{ for } \widetilde{\Phi}, \Phi \in  Act^{1}(F_2,\mathcal {S}^{2} )\}.$$

Let $\Phi$ be as above. It is clear that it is possible to take $\delta >0$ such that $\widetilde{\Phi}$ is a $C^{1}$ action with $d^{1}_S(\widetilde{\Phi}, \Phi)<\delta $ and:

\begin{enumerate}

 \item The action $\widetilde{\Phi}$ verifies the properties of definition \ref{definicion_delta}.
 \item The maps $\widetilde{\Phi}_a$ and $\widetilde{\Phi}_b$  are north-south pole  $C^{1}$ diffeomorphisms, $K_{\widetilde{\Phi}}$ is a Cantor set and for any $x,y\in K_{\widetilde{\Phi}}$, $x\neq y$, there exists $g\in F_2$ such that $d(\widetilde{\Phi}_g(x),\widetilde{\Phi}_g(y) )>\alpha$.\\
\end{enumerate}

By Lemma \ref{extension_lemma}, the action $\Phi$ is $C^{1}$-stable.

\vspace{.5cm}

\section{$C^{0}$-stability implies shadowing property.}
Here we will considerer the dynamical system $(F_2,M,\Phi)$ where $M$ is a compact manifold of dimension greater or equal two.

Recall that given $ g \in F_2 $ we denoted by $|g|$ the length of $g$.\\
Given $\delta >0$ and $n\in\N$, we say that $\{x_g\}_{|g|\leq n}$, $x_g\in M$, is a $\delta -n$ pseudotrajectory if $d(x_{sg},\Phi_s(x_g))<\delta$ for all $s\in\{ a,a^{-1},b,b^{-1}\}$, for all $g\in F_2 $ with $|g|\leq n$.

\begin{lemma}\label{shadowing_finito}
If for all $\varepsilon >0$ there exists $\delta >0$ such that if for all $n\in\N$  and for all  $\{x_g\}_{|g|\leq n}$ $\delta -n$, pseudotrajectory there exists $y\in M$ such that $d(\Phi_g(y), x_g )<\varepsilon $, $|g|\leq n$, then $\Phi$ has the shadowing property.
\end{lemma}

\pf

Given $\varepsilon$, choose $\delta$ as in the statement of the lemma. Let $\{x_g\}$ be a $\delta $ pseudotrajectory. For any $n\in \N$ we consider  $\{x_g\}_{|g|\leq n}$ , that is a $\delta -n$ pseudotrajectory included in  $\{x_g\}$. By hypotheses, there exist $y_n$ such that $d(\Phi_g(y_n), x_g )<\varepsilon $ for $|g|\leq n$.

Let $\{y_{n_{k}}\}$ be a subsequence of $\{y_n\}$ such that $y_{n_{k}}\to y$. Given $g\in F_2$ let $n_k$ be such that $|g|\leq n_k$ and  $d(\Phi_g(y), \Phi_g(y_{n_{k}})) <\varepsilon .$

  Then  $d(\Phi_g(y), x_g ) \leq  d(\Phi_g(y), \Phi_g(y_{n_{k}}))  + d(\Phi_g(y_{n_{k}}), x_g ) <2\varepsilon . $ And the proof is finished.\\

 For the proof of the next lemma see \cite[Lemma 13]{ns}.
\begin{lemma}\label{lema_difeo}
 Let M be a compact manifold
of dimension greater than or equal to two.
 Suppose a finite collection
$\{(p_i,q_i ) \in  M \times M \mbox {for }  i = 1,...,r\}$ is
specified together with a small $\lambda > 0$ such that
\begin{enumerate}
\item $ d(p_i,q_i) < \lambda$ for all $i$, and
\item if  $i\neq j$ then $p_i \neq  p_j$ and $q_i \neq q_j$.\\
Then there exists a diffeomorphism $f: M \to M$ such that\\
(a) $d(f,id) < 2\pi\lambda $, and\\
(b) $f(p_i) = q_i$   ($1 \leq  i \leq  r $).
\end{enumerate}

\end{lemma}

\begin{lemma}
If the action $\Phi $ is $C^{0}$-stable then the map $\Phi_{a^{-1}} \Phi_b$ is $C^{0}$-stable.
\end{lemma}

\pf
Let $f$ be $C^{0}$-close to $\Phi_{a^{-1}} \Phi_b $. So $\Phi_{a} f $ is $C^{0}$-close to $\Phi_b$. We considerer the action $\widetilde{\Phi} $ generated by $\Phi_{a} $ and $\Phi_{a} f $. As $\Phi$ is $C^{0}$-stable, there exists a semiconjugacy $h$ such that : $h\Phi_{a}=\Phi_{a}h$ and $h\Phi_{a} f = \Phi_{b}h$. Therefore $\Phi_{a}h f = \Phi_{b}h$ and $h f = \Phi_{a^{-1}}\Phi_{b}h$.\\

\begin{rk}\label{rk10}
As the map $\Phi_{a^{-1}} \Phi_b$ is $C^{0}$-stable there is no an open set $U\subset M$ such that every point of $U$ is a periodic point, with all of them with the same period.
\end{rk}

\begin{lemma}\label{lema4}

 Let $\Phi$ be a $C^{0}$-stable action. Given $\eta >0$  there exists $\delta_0>0$ such that if  $\{  x_g  \}$ is a $\delta$-pseudotrajectory for $\Phi$ with $\delta <\delta_0$, then there exists  $\{ \widetilde{x}_g \} $
 such that:
\begin{itemize}
\item $d(x_g, \widetilde{x}_g)<\eta $ and

 \item $\Phi_a( \widetilde{x}_{a^{-1}g} )=     \Phi_b(\widetilde{x}_{b^{-1}g} )$  for all $g\in F_2.$
\end{itemize}

\end{lemma}

\pf

Let $\{  x_g \} $ be a $\delta$-pseudotrajectory.
We will define an equivalence relation in  $\{  x_g \} $. We say that $  x_g $ is equivalent to $  x_{g '} $ if there exists  $n \in \Z$ such that $ g= (b^{-1}a)^n g' $. Let $ [ x_g]   $ be the class of element of $x_g$. Note that there exists $g_0 \in F_2$ with minimum length such that  $  [ x_g]   = [ x_{g_0}]   =  \{       x_{(b^{-1}a)^ng_0}   : n\in Z       \}$. Therefore  the class $  [ x_g]  $ can be thought as a sequence $\{ z_n \}$ such that $ z_n = x_{(b^{-1}a)^ng_0} $, $n\in \Z$.

Let's prove that $\{ z_n \}$ is a $\delta +\delta_1  $-pseudotrajectory for the map $\Phi_{b^{-1}} \Phi_a$, where $\delta_1 (\delta)=Max\{ d(\Phi_{b^{-1}}(x) ,\Phi_{b^{-1}}(y) ):\  \ d(x,y)<\delta     \}$.
  $$ d(   \Phi_{b^{-1}} \Phi_a (z_n), z_{n+1}  )  =d(  \Phi_{b^{-1}} \Phi_a (x_{(b^{-1}a)^ng_0}), x_{(b^{-1}a)^{n+1}g_0}     )\leq $$
         $$     d(   \Phi_{b^{-1}} \Phi_a (   x_{(b^{-1}a)^ng_0}      ),  \Phi_{b^{-1}}( x_{a(b^{-1}a)^n g_0}     )  )      +d(     \Phi_{b^{-1}}( x_{a(b^{-1}a)^n g_0}     )     ,           x_{(b^{-1}a)^{n+1}g_0}                                       )            < \delta_1 (\delta) + \delta  .   $$

As the map $\Phi_{b^{-1}} \Phi_a $ is $C^{0}$-stable, then it has the shadowing property  (see \cite{w} ). So given $\eta>0$ there exists $\delta^{'}_0$ such that if  $\delta_1(\delta) + \delta <\delta^{'}_0$ then there exists $y$ such that\\ $d(( \Phi_{b^{-1}} \Phi_a)^n (y), z_n     ) <\eta    $.

Therefore, for each $[ x_{g_{0}}] $ there exists $y_{g_{0}}\in M$ such that\\ $d(( \Phi_{b^{-1}} \Phi_a)^n (y_{g_0}),      x_{(b^{-1}a)^n g_0}     ) <\eta    $.

Now, we define $\{  \widetilde{x}_g  \} $  such that $ \widetilde{x}_g = ( \Phi_{b^{-1}} \Phi_a)^n (y_{g_0})$  if $g=(b^{-1}a)^ng_0$.

Note that $x_{a^{-1}g} $ is equivalent to $  x_{b^{-1}g} $ , so $\Phi_{b^{-1}a }( \widetilde{x}_{a^{-1}g}  )= \widetilde{x}_{b^{-1}g}   $. Therefore $\Phi_{a }( \widetilde{x}_{a^{-1}g}  )=\Phi_{b } (\widetilde{x}_{b^{-1}g})$ and $d(x_g, \widetilde{x}_g  )<\eta$. So given $\eta >0$  it is enough to take $\delta_0$ such that $\delta_1(\delta_0)+\delta_0<\delta^{'}_0$.\\

{\bf{Proof of Theorem A}}.\\

  Given $\eta >0$, let $\delta_0$ be given by Lemma \ref{lema4}. Let $\{x_g\}$ be a $\delta$  pseudotrajectory with $\delta <\delta_0$.  For each $n\in \N$, let   $\{x_g\}_{|g|\leq n}$  be, a $\delta -n$  pseudotrajectory. Again, by Lemma \ref{lema4} there exists $\{ \widetilde{x}_g\}$ such that   $d(x_g, \widetilde{x}_g)<\eta $ for $|g|\leq n$. By construction, $ \widetilde{x}_g = ( \Phi_{b^{-1}} \Phi_a)^m (y_{g_{0}})$. As $|g|\leq n$, by Remark \ref{rk10}, it is possible to take $y_{g_0}$ such that $\widetilde{x}_g \neq \widetilde{x}_{g^{'}}$ if $g\neq g^{'}$ with $|g|\leq n$ and $|g^{'}|\leq n$.\\
By Lemma \ref{shadowing_finito} it is enough to prove the Theorem A for a $\delta -n$  pseudotrajectory.

Given $\varepsilon >0$ let

\begin{itemize}

\item $\rho >0$, $\rho<\varepsilon/2$, such that, if $\widetilde{\Phi}$ is an action with $d^{0}_S( \widetilde{\Phi},\Phi  )<\rho$ then there exists $h:M\to M$, $d(h,id)<\varepsilon /2$ and $h\widetilde{\Phi}_g=\Phi_g h$ for all $g\in F_2$.

\item $\eta >0$, $\eta<\varepsilon/2$. We take $\delta <\delta_0$ small enough such that for any  $\{\widetilde{x}_g\}$  as above, ($|g|\leq n$) with $d(x_g, \widetilde{x}_g)<\eta $ there exists a diffeomorphisms $f:M\to M$, given by Lemma \ref{lema_difeo}, such that:

\begin{enumerate}

\item\noindent\begin{equation}\label{pp}
 \noindent    \  f(\Phi_{a }( \widetilde{x}_{a^{-1}g}  ))=f(\Phi_{b } (\widetilde{x}_{b^{-1}g}))=\widetilde{x}_g,
\end{equation}
\item If $\widetilde{\Phi}$ is the action generated by $\widetilde{\Phi}_a=f\Phi_a$ and  $\widetilde{\Phi}_b=f\Phi_b$, then $$ \ d^{0}_S(\widetilde{\Phi}, \Phi)<\rho . $$

\end{enumerate}

\end{itemize}

  Let's prove that $\widetilde{\Phi}_s   ( \widetilde{x}_g    )= \widetilde{x}_{sg}$ for all $s\in \{a,a^{-1},b,b^{-1} \}$, for all $g\in F_2$. It is enough to prove that  $\widetilde{\Phi}_s   ( \widetilde{x}_g    )= \widetilde{x}_{sg}$ for all $s\in \{a,b \}$, for all $g\in F_2$.
Suppose $s=a$. If $s=b$ the proof is analogous.
$$\widetilde{\Phi}_a   ( \widetilde{x}_g    )= \widetilde{\Phi}_a   ( \widetilde{x}_{a^{-1}ag}    )               \stackrel{\mbox{\scriptsize
\begin{tabular}{c}
$(\ref{pp})$
\end{tabular}}}   =\widetilde{x}_{ag}  . $$

            As $d^{0}_S( \widetilde{\Phi},\Phi  )<\rho$, then there exists $h:M\to M$ such that $h\widetilde{\Phi}_g =\Phi_g h$ for all $g\in F_2$, with $d(h(x),x)<\varepsilon /2$ . So

$$  d(\Phi_g( h(\widetilde{x}_e)), x_g     ) = d(h\widetilde{\Phi}_g(\widetilde{x}_e), x_g     ) = d(h(\widetilde{x}_g), x_g     )$$

$$\leq    d(h(\widetilde{x}_g), \widetilde{x}_g     ) + d(\widetilde{x}_g, x_g     )              <   \varepsilon/2+\eta  <\varepsilon .$$

\end{document}